\documentclass[journal]{IEEEtran}
\usepackage{amsfonts}
\usepackage{amsthm,amsmath, amssymb}
\usepackage{graphicx}
\usepackage{verbatim}
\usepackage{bbm}
\usepackage{rotating}

\usepackage{algorithmic}
\usepackage{algorithm}

\usepackage{tikz}
\usetikzlibrary{matrix}

\newcommand{\bo}[1]{\boldsymbol{#1}}

\def\cov{\mathrm{Cov}}
\newcommand{\E}{\mbox{E}}
\newcommand{\var}{\mbox{Var}}
\newtheorem{theorem}{Theorem}

\newtheorem{corollary}{Corollary}

\newtheorem{assumption}{Assumption}
\newtheorem{remark}{Remark}
\newcommand{\bmat}{\begin{pmatrix}}
\newcommand{\emat}{\end{pmatrix}}
\usepackage{color}

\def\vec{\mathrm{vec}}

\graphicspath{ {./Figures/} }

\newcommand\numberthis{\addtocounter{equation}{1}\tag{\theequation}}

\begin{document}
%
\title{JADE for Tensor-Valued Observations}
%
%
%

\author{Joni~Virta, Bing~Li, Klaus~Nordhausen and Hannu~Oja
\thanks{J. Virta, K. Nordhausen and H. Oja are with the Department of Mathematics and Statistics, University of
Turku, 20014 Turku, Finland
(e-mail: joni.virta@utu.fi).}
\thanks{B. Li is with the Department of Statistics, Pennsylvania State University, 326 Thomas Building, University Park, Pennsylvania 16802, USA.}
}

%
%

\markboth{}%
{Virta \MakeLowercase{\textit{et al.}}: Tensor observation JADE}
%



\maketitle

\begin{abstract}
Independent component analysis is a standard tool in modern data analysis and numerous different techniques for applying it exist. The standard methods however quickly lose their effectiveness when the data are made up of structures of higher order than vectors, namely matrices or tensors (for example, images or videos), being unable to handle the high amounts of noise. Recently, an extension of the classic fourth order blind identification (FOBI) specifically suited for tensor-valued observations was proposed and showed to outperform its vector version for tensor data. In this paper we extend another popular independent component analysis method, the joint approximate diagonalization of eigen-matrices (JADE), for tensor observations. In addition to the theoretical background we also provide the asymptotic properties of the proposed estimator and use both simulations and real data to show its usefulness and superiority over its competitors.

\end{abstract}

\begin{IEEEkeywords}
Independent component analysis, multilinear algebra, kurtosis, limiting normality, minimum distance index.
\end{IEEEkeywords}


%
\IEEEpeerreviewmaketitle

\section{Introduction}\label{sec:intro}

The following presentation relies on multilinear algebra and before the actual ideas can be described we first review some key properties of tensors and matrices needed later.

A tensor of $r$th order $\textbf{X} \in \mathbb{R}^{p_1 \times \cdots \times p_r}$ can be seen as a higher order analogy of vectors and matrices. Whereas a matrix can be viewed either as a collection of rows or that of columns, a tensor of $r$th order has in total $r$ \textit{modes}. The \textit{$m$-mode vectors} of a tensor are given by letting the $m$th index vary while keeping all other indices fixed, $m = 1,\ldots,r$. A tensor $\textbf{X} \in \mathbb{R}^{p_1 \times \cdots \times p_r}$ thus contains $\rho_m := \Pi_{s \neq m}^r p_s$ $m$-mode vectors of length $p_m$. The opposite construct, fixing a single index $i_m$ and varying the others, then gives what we call the \textit{$m$-mode faces} of a tensor. The number of $m$-mode faces then totals $p_m$ and each is a tensor of size $p_1 \times \cdots \times p_{m-1} \times p_{m+1} \times \cdots \times p_r$.

For representing tensor contraction, or summation, we use the Einstein summation convention in which a twice-appearing index in a product implies summation over the range of the index. 
For example, for a tensor $\textbf{X} = \{x_{i_1 i_2 i_3}\}$ we have
\[x_{i_1 i_2 j} x_{i_1 i_2 k} := \sum_{i_1=1}^{p_1} \sum_{i_2=1}^{p_2} x_{i_1 i_2 j} x_{i_1 i_2 k}.\]
Two special cases of tensor contraction prove especially useful for us. The product $\textbf{X} \odot_m \textbf{A}$ of tensor $\textbf{X} \in \mathbb{R}^{p_1 \times \cdots \times p_r}$ with a matrix $\textbf{A} \in \mathbb{R}^{p_m \times p_m}$, $m=1,\ldots,r$, is defined as the ${p_1 \times \cdots \times p_r}$-dimensional tensor with the elements
\begin{equation}\label{eq:tensorbymatrix}
\left(\textbf{X} \odot_m \textbf{A}\right)_{i_1 \ldots i_r} = x_{i_1 \ldots i_{m-1} j_m i_{m+1} \ldots i_r} a_{i_m j_m}.
\end{equation}
That is, the multiplication $\textbf{X} \odot_m \textbf{A}$ linearly transforms $\textbf{X}$ from the direction of the $m$th mode without changing the size of the tensor. The operation can alternatively be viewed as applying the linear transformation given by $\textbf{A}$ separately to each $m$-mode vector of the tensor. The second useful product, $\textbf{X} \odot_{-m} \textbf{Y}$, of two tensors of the same size, $\textbf{X}, \textbf{Y} \in \mathbb{R}^{p_1 \times \cdots \times p_r}$ is defined as the $p_m \times p_m$-dimensional matrix with the elements
\begin{equation}\label{eq:tensorbytensor}
\left(\textbf{X} \odot_{-m} \textbf{Y}\right)_{jk} = x_{i_1 \ldots i_{m-1} j i_{m+1} \ldots i_r} y_{i_1 \ldots i_{m-1} k i_{m+1} \ldots i_r}.
\end{equation}
The special case $\textbf{X} \odot_{-m} \textbf{X}$ provides higher order counterparts for the products of a vector $\textbf{x} \in \mathbb{R}^{p_1}$ or a matrix $\textbf{X} \in \mathbb{R}^{p_1 \times p_2}$ with itself, such as $\textbf{x}\textbf{x}^T$, $\textbf{X} \textbf{X}^T$ or $\textbf{X}^T \textbf{X}$, and proves useful in defining the ``covariance matrix'' of a tensor.

Finally, define the vectorization $\vec (\textbf{X}) \in \mathbb{R}^{p_1 \cdots p_r}$ of a tensor $\textbf{X} \in \mathbb{R}^{p_1 \times \cdots \times p_r}$ as the stacking of the elements $x_{i_1 \ldots i_r}$ in such a way that the leftmost index goes through its cycle the quickest and the rightmost index the slowest. Then it holds for a tensor $\textbf{X} \in \mathbb{R}^{p_1 \times \cdots \times p_r}$ and matrices $\bo{A}_1 \in \mathbb{R}^{p_1 \times p_1},\ldots,\bo{A}_r \in \mathbb{R}^{p_r \times p_r}$ that
\[\vec (\textbf{X} \odot_1 \bo{A}_1 \cdots \odot_r \bo{A}_r) = (\textbf{A}_r \otimes \cdots \otimes \textbf{A}_1) \vec (\textbf{X}),\]
where $\otimes$ is the Kronecker product.

In this paper we assume that the tensor-valued i.i.d. random elements $\textbf{X}_i \in \mathbb{R}^{p_1 \times \cdots \times p_r}$, $i=1,\ldots,n$, are observed from the recently suggested \cite{virta2015mfobi} \textit{tensor independent component model}:
\begin{align}\label{eq:tica_model}
\textbf{X} = \bo{\mu} + \textbf{Z} \odot_1 \bo{\Omega}_1 \cdots \odot_r \bo{\Omega}_r,
\end{align}
where $\bo{\Omega}_1 \in \mathbb{R}^{p_1 \times p_1},\ldots,\bo{\Omega}_r \in \mathbb{R}^{p_r \times p_r}$ are full rank \textit{mixing matrices}, $\bo{\mu} \in \mathbb{R}^{p_1 \times \cdots \times p_r}$ is the location center, and $\textbf{Z} \in \mathbb{R}^{p_1 \times \cdots \times p_r}$ is an unobserved random tensor. The model \eqref{eq:tica_model} is further equipped with the following assumptions.
\begin{assumption}\label{assu:ind}
The components of $\textbf{Z}$ are mutually independent.
\end{assumption}
\begin{assumption}\label{assu:stand}
The components of $\textbf{Z}$ are standardized in the sense that $E[\vec (\textbf{Z})] = \textbf{0}$ and $\cov [ \vec (\textbf{Z})] = \textbf{\em I}$.
\end{assumption}
\begin{assumption}\label{assu:gauss}
For each $m = 1,\ldots,r$, at most one $m$-mode face of $\textbf{Z}$ consists entirely of Gaussian components.
\end{assumption}

Assumption \ref{assu:stand} implies that
$E[\textbf{X}]=\bo\mu$ and that
\[ \cov [ \vec (\textbf{X})] = (\bo{\Omega}_r \bo{\Omega}_r^T) \otimes \cdots \otimes  (\bo{\Omega}_1 \bo{\Omega}_1^T) \]
has the so-called Kronecker structure. Assumption \ref{assu:gauss} is a tensor analogy for the usual vector independent component model assumption on maximally one Gaussian component and without it some column blocks of some of the matrices $\bo{\Omega}_1,\ldots,\bo{\Omega}_{r}$ could be identifiable only up to a rotation. After the above assumptions we can still freely change the signs and orders of the columns of all $\bo{\Omega}_1,\ldots,\bo{\Omega}_{r}$, or multiply any $\bo{\Omega}_s$ by a constant and divide any $\bo{\Omega}_t$ by the same constant, but this indeterminacy is acceptable in practice. The model along with its assumptions now provides a natural extension for the standard independent component model which is obtained as a special case when $r = 1$.

Alternatively, the model can be seen as an extension of the general location-scatter model for tensor-valued data, which is equivalent to \eqref{eq:tica_model} with only Assumption \ref{assu:stand} and is often, for $r=1,2$, combined with the assumption on Gaussianity or sphericity of $\vec (\textbf{Z})$. Under the location-scatter model the covariance matrix of $\vec(\textbf{X})$ again has the above Kronecker structure. In addition to requiring less parameters to estimate than a full $p_1 \cdots p_r \times p_1 \cdots p_r$ covariance matrix, the assumption on Kronecker structure is a natural choice in many applications, see e.g. \cite{werner2008estimation}. For the estimation of covariance parameters under the assumption on Kronecker structure in the matrix case, $r=2$, see \cite{srivastava2008models, wiesel2012geodesic, sun2015robust}. For the general tensor Gaussian distribution and the estimation of its parameters see \cite{ohlson2013multilinear, manceur2013maximum}.

The extension of dimension reduction methods from vector to matrix or tensor observations is in signal processing usually
approached via different tensor decompositions such as the CP-decomposition and the Tucker decomposition. A review of them along with a plethora of references for applications is given in \cite{kolda2009tensor}, see also \cite{lu2011survey} for more applications. For examples of particular dimension reduction methods incorporating matrix or tensor predictors, see e.g. \cite{vasilescu2005multilinear, zhang2008directional, virta2015mfobi} for independent component analysis, \cite{li2010dimension, pfeiffer2012sufficient, xue2014sufficient, ding2015tensor} for sufficient dimension reduction and \cite{ding2014dimension, greenewald2014robust} for principal components analysis-based techniques. More references are also given in \cite{li2010dimension, virta2015mfobi}.


In tensor independent component analysis the objective is to estimate, based on the sample $\textbf{X}_1, \ldots, \textbf{X}_n$, some unmixing matrices $\bo{\Phi}_1,\ldots,\bo{\Phi}_r$ such that $\textbf{X} \odot_1 \bo{\Phi}_1 \cdots \odot_r \bo{\Phi}_r$ has mutually independent components. A na{\"i}ve method for accomplishing this would be to vectorize the observations and resort to some standard method of independent component analysis, but in doing so the resulting estimate lacks the desired Kronecker structure. In addition, vectorizing and using standard tools meant for vector-valued data requires the stronger, component-wise version of Assumption~\ref{assu:gauss}, inflates the number of parameters and can make the dimension of the data too large for standard methods to handle. To circumvent this, \cite{vasilescu2005multilinear, zhang2008directional, virta2015mfobi} proposed estimating an unmixing matrix separately for each of the modes and \cite{virta2015mfobi} presented an extension of the classic fourth order blind identification (FOBI) \cite{cardoso1989source} for tensor observations called TFOBI.

In the vector independent component model, $\textbf{x}=\bo \mu+\bo\Omega \textbf{z}$, the standardized vector $\textbf{x}_{st} := \cov [\textbf{x}]^{-1/2} (\textbf{x} - E[\textbf{x}])$ equals $\textbf{U}\textbf{z}$ for some orthogonal matrix $\textbf{U}$, see e.g. \cite{miettinen2014fourth}. In FOBI  the rotation $\textbf{U}$ is then found using the eigendecomposition of the matrix of fourth moments $\textbf{B} := E[\textbf{x}_{st} \textbf{x}_{st}^T \textbf{x}_{st} \textbf{x}_{st}^T]$. This same approach is taken in TFOBI by performing both steps of the procedure, the standardization and the rotation, on all $r$ modes of $\textbf{X}$. Assuming centered $\textbf{X}$, in \cite{virta2015mfobi} the $m$-mode covariance matrices,
\begin{align}\label{eq:mjade_covs}
\bo{\Sigma}_m(\textbf{X}) := \rho_m^{-1} E\left[\textbf{X} \odot_{-m} \textbf{X} \right], \quad m = 1, \ldots, r,
\end{align}
are first used to standardize the observations as $\textbf{X}_{st} := \textbf{X} \odot_1 \bo{\Sigma}{}_1^{-1/2} \cdots \odot_r \bo{\Sigma}{}_r^{-1/2}$. The tensor $\textbf{Z}$ is then found by rotating $\textbf{X}_{st}$ from all $r$ modes and the rotation matrices can be found from the eigendecompositions of the $m$-mode matrices of fourth moments:
\[\textbf{B}_m := \rho_m^{-1} E\left[(\textbf{X}_{st} \odot_{-m} \textbf{X}_{st}) (\textbf{X}_{st} \odot_{-m} \textbf{X}_{st}) \right]. \]

Another widely used independent component analysis method for vector-valued data, called the joint approximate diagonalization of eigen-matrices (JADE) \cite{cardoso1993blind}, also uses fourth moments to estimate the required final rotation but utilizes them in the form of cumulant matrices (assuming $E(\textbf{x})=\textbf{0}$)
\begin{align*}
\textbf{C}^{ij}(\textbf{x}) &:= E\left[x_i x_j \cdot \textbf{x} \textbf{x}^T\right] -  E[x_i x_j] E\left[\textbf{x} \textbf{x}^T\right]  \numberthis \label{eq:jade_cij 1} \\
&-  E\left[x_i \cdot \textbf{x} \right]  E\left[x_j \cdot \textbf{x}^T \right] -  E\left[x_j \cdot \textbf{x} \right]  E\left[x_i \cdot \textbf{x}^T \right].
\end{align*}
The final rotation from $\textbf{x}_{st}$ to $\textbf{z}$ is in JADE obtained by jointly diagonalizing the matrices
\begin{align}\label{eq:jade_cij}
\textbf{C}^{ij}(\textbf{x}_{st}) = E\left[x_{st,i} x_{st,j} \cdot \textbf{x}_{st} \textbf{x}_{st}^T\right] - \delta_{ij} \textbf{I} - \textbf{E}^{ij} - \textbf{E}^{ji},
\end{align}
where $\textbf{E}^{ij}$ is a matrix with a single one as element $(i, j)$ and zeroes elsewhere and $\delta_{ij}$ is the Kronecker delta. Compared to FOBI which only uses $p(p+1)/2$ sums of fourth joint moments of $\textbf{x}_{st}$ JADE thus has a clear advantage in using all possible fourth joint cumulants of $\textbf{x}_{st}$ in the estimation of the rotation matrix.

Because of the well-known fact that JADE outperforms FOBI in most cases it is natural to expect that the extension of JADE to tensor-valued data would similarly be superior to TFOBI. This is indeed the case, and in the following sections we formulate the tensor joint diagonalization of eigen-matrices (TJADE) which is obtained from JADE by applying very much the same extensions as required when moving from FOBI to TFOBI. We first briefly discuss the standard vector-valued independent component model and review the theory and assumptions behind the original JADE in Section~\ref{sec:jade}. The corresponding aspects of TJADE are presented in Section~\ref{sec:mjade} and the asymptotical properties of both methods in Section~\ref{sec:asymp}. Simulations comparing TJADE to TFOBI and both the original JADE and original FOBI are presented in Section~\ref{sec:simu} along with a real data example and we close in Section~\ref{sec:disc} with some discussion. The proofs can be found in Appendix~\ref{sec:appen}.

\section{Original JADE}\label{sec:jade}

The original JADE assumes that the vector-valued observations are generated by the vector independent component model
\begin{align}\label{eq:vicm_model}
\textbf{x}_i = \bo{\mu} + \bo{\Omega} \textbf{z}_i, \quad i=1,\ldots,n,
\end{align}
where the mixing matrix $\bo{\Omega} \in \mathbb{R}^{p \times p}$ has full rank, $\bo{\mu} \in \mathbb{R}^p$ and the i.i.d. random vectors $\textbf{z}_i \in \mathbb{R}^p$ have mutually independent components standardized to have zero means and unit variances. To ensure the existence of the JADE solution we have to further assume that at most one of the independent components has zero excess kurtosis \cite{miettinen2014fourth}.

Assuming next that the data are centered, that is, $E[\textbf{x}] = \textbf{0}$, we standardize the vectors as $\textbf{x}_{st} = \bo{\Sigma}^{-1/2} \textbf{x}$. The standardized vectors can be shown to satisfy $\textbf{x}_{st} = \textbf{U}\textbf{z}$ for some orthogonal matrix $\textbf{U}$, see for example \cite{miettinen2014fourth}. To estimate $\textbf{U}$, JADE uses the cumulant matrices $\textbf{C}^{ij}(\textbf{x}_{st})$, $i,j=1,\ldots,p$, in \eqref{eq:jade_cij}. Under the independent component model the cumulant matrices can be shown to satisfy, for all $i,j=1,\ldots,p$,
\begin{align}\label{eq:jade_cumulant}
\textbf{C}^{ij}(\textbf{x}_{st}) = \textbf{U}  \left( \sum_{k=1}^p u_{ik} u_{jk} \kappa_k \textbf{E}^{kk} \right) \textbf{U}^T,
\end{align}
where $\kappa_k := E(z_k^4) - 3$, the excess kurtosis of the $k$th component, and $u_{ab}$ are the components of $\textbf{U}$. The expression in \eqref{eq:jade_cumulant} is the eigendecomposition of $\textbf{C}^{ij}(\textbf{x}_{st})$ and thus any single matrix    $\textbf{C}^{ij}(\textbf{x}_{st})$ could be used to find $\textbf{U}$. However, to use all the information available in the fourth joint cumulants, JADE simultaneously (approximately) diagonalizes them all, that is, finds $\textbf{U}^T$ as
\begin{align}\label{eq:jade_simul}
\textbf{U}^T = \underset{\textbf{U}: \ \textbf{U}^T\textbf{U} = \textbf{I}}{\mbox{argmax}} \sum_{i=1}^p \sum_{j=1}^p \| \mbox{diag}(\textbf{U} \textbf{C}^{ij}(\textbf{x}_{st}) \textbf{U}^T) \|^2.
\end{align}
Optimization problems of type \eqref{eq:jade_simul} are so-called joint diagonalization problems for which many algorithms exist, see \cite{belouchrani1997blind} for discussion and one particular algorithm.

In \cite{miettinen2014fourth}, a thorough analysis of the statistical properties of JADE is given and it is shown there that the JADE estimator is an \textit{independent component functional}, that is, the resulting components are invariant up to sign-change and permutation under affine transformations to the original data.

\section{Tensor joint approximate diagonalization of eigen-matrices}\label{sec:mjade}

In formulating TJADE we assume that the data are generated by the tensor independent component model \eqref{eq:tica_model} and satisfy Assumptions \ref{assu:ind}, \ref{assu:stand} and \ref{assu:gauss}. Assuming $E[\textbf{X}]= \textbf{0}$, we next go separately through the tensor analogies of the standardization and rotation steps of the original JADE.

\subsection{Standardization step}

We take the same approach for standardization of $\textbf{X}$ as in \cite{virta2015mfobi}, that is, use the $m$-mode covariance matrices, $\bo{\Sigma}_1,\ldots,\bo{\Sigma}_r$, to standardize $\textbf{X}$ simultaneously from all $r$ modes. This gives us the standardized tensor
\[ \textbf{X}_{st} := \textbf{X} \odot_1 \bo{\Sigma}{}_1^{-1/2} \cdots \odot_r \bo{\Sigma}{}_r^{-1/2}. \]
where, for the asymptotics, we assume that the standardization functionals $\bo{\Sigma}{}_m^{-1/2}$, $m=1,\ldots,r$, are chosen to be symmetric, see e.g. \cite{ilmonen2012invariant}. Estimates $\hat{\bo{\Sigma}}_1,\ldots,\hat{\bo{\Sigma}}_r$ of the $m$-mode covariance matrices are obtained by applying \eqref{eq:mjade_covs} to the empirical distribution of $\textbf{X}$. 
The next step towards $\textbf{Z}$ is guided by Theorem 5.3.1 in \cite{virta2015mfobi} which states that
\begin{align} \label{eq:mjade_stand}
\textbf{X}_{st} = \tau \cdot \textbf{Z} \odot_1 \textbf{U}_1 \cdots \odot_r \textbf{U}_r,
\end{align}
for some orthogonal matrices $\textbf{U}_1 \in \mathbb{R}^{p_1 \times p_1},\ldots,\textbf{U}_r \in \mathbb{R}^{p_r \times p_r}$ and for $\tau = (\prod_{i=1}^m p_m^{1/2})^{r-1} \|\bo\Omega_r\otimes \cdots \otimes \bo\Omega_1\|_F^{1-r}$, where $\| \cdot \|_F$ is the Frobenius norm.

%
%

\subsection{Rotation step}

We extend the cumulant matrices by noting that the operation $\odot_{-m}$ provides an $m$-mode analogy for the outer product of vectors. By writing the random quantity $x_i x_j \cdot \textbf{x} \textbf{x}^T$ in \eqref{eq:jade_cij 1} with outer products either as $\textbf{e}_i^T \textbf{x} \textbf{x}^T \textbf{e}_j \cdot \textbf{x} \textbf{x}^T$ or as $\textbf{x} \textbf{x}^T \textbf{e}_i \textbf{e}_j^T\textbf{x} \textbf{x}^T$ two straightforward tensor $m$-mode analogies for the matrix of fourth cumulants $\textbf{C}^{ij}$, $i,j=1,\ldots,p_m$, in \eqref{eq:jade_cij 1} are then given by
\begin{align}
\begin{split}
\textbf{C}_{1,m}^{ij}(\textbf{X}) &= \rho_m^{-1} E\left[ \textbf{e}_i^T(\textbf{X} \odot_{-m} \textbf{X}) \textbf{e}_j \cdot (\textbf{X} \odot_{-m} \textbf{X}) \right] \\
&- \rho_m^{-1} E\left[ \textbf{e}_i^T(\textbf{X}^* \odot_{-m} \textbf{X}^*) \textbf{e}_j \cdot (\textbf{X} \odot_{-m} \textbf{X}) \right] \\
&- \rho_m^{-1} E\left[ \textbf{e}_i^T(\textbf{X}^* \odot_{-m} \textbf{X}) \textbf{e}_j \cdot (\textbf{X}^* \odot_{-m} \textbf{X}) \right]  \\
&- \rho_m^{-1} E\left[ \textbf{e}_i^T(\textbf{X}^* \odot_{-m} \textbf{X}) \textbf{e}_j \cdot (\textbf{X} \odot_{-m} \textbf{X}^*) \right],
\end{split}  \label{eq:mjade_cij1}
\end{align}
and
\begin{align}
\begin{split}
\textbf{C}_{2,m}^{ij}(\textbf{X}) &= \rho_m^{-1} E\left[ (\textbf{X} \odot_{-m} \textbf{X}) \textbf{E}^{ij} (\textbf{X} \odot_{-m} \textbf{X}) \right] \\ &- \rho_m^{-1} E\left[ (\textbf{X}^* \odot_{-m} \textbf{X}^*) \textbf{E}^{ij} (\textbf{X} \odot_{-m} \textbf{X}) \right]  \\&- \rho_m^{-1} E\left[ (\textbf{X}^* \odot_{-m} \textbf{X}) \textbf{E}^{ij} (\textbf{X}^* \odot_{-m} \textbf{X}) \right]  \\
&- \rho_m^{-1} E\left[ (\textbf{X}^* \odot_{-m} \textbf{X}) \textbf{E}^{ij} (\textbf{X} \odot_{-m} \textbf{X}^*) \right],
\end{split}\label{eq:mjade_cij2}
\end{align}
with $m = 1,\ldots,r$, where $\textbf{X}^*$ is an independent copy of $\textbf{X}$. Theoretically, a third way to generalize the idea is obtained by considering $\textbf{x} \textbf{x}^T \textbf{e}_j \textbf{e}_i^T\textbf{x} \textbf{x}^T$. However, that would be redundant as the resulting set of matrices for $i,j=1,\ldots,p_m$ is the same as with \eqref{eq:mjade_cij2} and the individual matrices can be obtained by just reversing $i$ and $j$ in \eqref{eq:mjade_cij2}. Naturally, for vector observations, $r=1$, both \eqref{eq:mjade_cij1} and \eqref{eq:mjade_cij2} are equivalent.

\def\lo#1{_{#1}}
\def\hi#1{^{#1}}
\def\tbf{\textbf}

Define next for the model \eqref{eq:tica_model} its kurtosis tensor $\boldsymbol{\kappa} \in \mathbb{R}^{p_1 \times \cdots \times p_r}$ as $(\boldsymbol{\kappa})_{i_1 \ldots i_r} := E[z_{i_1 \ldots i_r}^4] - 3$ and its $m$-mode average kurtosis vector as $\bar{\bo{\kappa}}^{(m)} := (\bar \kappa{}_ 1^{(m)},\ldots,\bar \kappa{}_{p_m}^{(m)})$, where $\bar \kappa _ k ^{(m)}$ is the average of the excess kurtoses of the random variables in the $k$th $m$-mode face of the tensor $\textbf{Z}$, $k=1,\ldots,p_m$. The following theorem then shows that \eqref{eq:mjade_cij1} and \eqref{eq:mjade_cij2} actually serve in TJADE the same purpose as their vector counterpart in JADE.

\begin{theorem}\label{theo:cumulant_diag} If $\tau$, $\textbf{U}_1$, \ldots, $\textbf{U}_r$ are as defined in (\ref{eq:mjade_stand}), then, for  $c=1,2$ and $m = 1,\ldots,r$, the matrices of fourth cumulants $\textbf{C}_{c,m}^{ij}, i, j = 1,\ldots,p$ satisfy 
\[\textbf{C}_{c, m}^{ij}(\textbf{X}_{st}) = \tau^4\cdot \textbf{U}_m \left( \sum_{k=1}^{p_m} u^{(m)}_{ik} u^{(m)}_{jk} \bar{\kappa}{}^{(m)}_k \textbf{E}^{kk} \right) \textbf{U}_m^T. \]
\end{theorem}


According to Theorem~\ref{theo:cumulant_diag}, $\textbf{U}_m^T$ simultaneously diagonalizes all matrices $\textbf{C}_{c, m}^{ij}(\textbf{X}_{st})$, $i,j=1,\ldots,p_m$, regardless of $c$, giving two straightforward ways of estimating the $m$-mode rotation $\textbf{U}_m$ using \eqref{eq:jade_simul} with $\textbf{C}^{ij}(\textbf{x}_{st})$ replaced by $\textbf{C}_{c, m}^{ij}(\textbf{X}_{st})$ for the chosen value of $c$. However, in estimating an individual matrix $\textbf{C}_{c, m}^{ij}(\textbf{X}_{st})$ in \eqref{eq:mjade_cij1} or \eqref{eq:mjade_cij2} we have to estimate four matrices in total, the last two of which are costly to estimate because of the independent copies $\textbf{X}^*$. Using the method of the proof of Theorem \ref{theo:cumulant_diag} one can show that, analogously to the vector-valued case,
\begin{align*}
\textbf{C}_{1, m}^{ij}(\textbf{X}_{st}) &= \textbf{B}_{1, m}^{ij} -
\bo{\Xi}_m \left(\delta_{ij} \rho_m \textbf{I} + \textbf{E}^{ij} + \textbf{E}^{ji}\right) \bo{\Xi}_m^T,
\end{align*}
where $\textbf{B}_{1, m}^{ij} := \rho_m^{-1} E\left[ \textbf{e}_i^T (\textbf{X}_{st} \odot_{-m} \textbf{X}_{st}) \textbf{e}^j \cdot (\textbf{X}_{st} \odot_{-m} \textbf{X}_{st}) \right]$  and $\bo{\Xi}_m := \rho_m^{-1} E\left[ \textbf{X}_{st} \odot_{-m} \textbf{X}_{st} \right] = \tau^2 \textbf{I}$, which provides a natural estimator for $\tau^2$. Similarly
\begin{align*}
\textbf{C}_{2, m}^{ij}(\textbf{X}_{st}) &= \textbf{B}_{2, m}^{ij} -
\bo{\Xi}_m \left(\delta_{ij} \textbf{I} + \rho_m \textbf{E}^{ij} + \textbf{E}^{ji}\right) \bo{\Xi}_m^T,
\end{align*}
where $\textbf{B}_{2, m}^{ij} := \rho_m^{-1} E\left[ (\textbf{X}_{st} \odot_{-m} \textbf{X}_{st}) \textbf{E}^{ij} (\textbf{X}_{st} \odot_{-m} \textbf{X}_{st}) \right]$ and $\bo{\Xi}_m$ is as above.

\def\ali{&\,}
\def\inv{\hi {-1}}

Natural estimates for the previous matrices are provided by
\begin{align}\label{eq:mjade_kumulant1}
\hat{\textbf{C}}{}_{1, m}^{ij} := \hat{\textbf{B}}{}_{1, m}^{ij} - \hat{\bo{\Xi}}_m \left(\delta_{ij} \rho_m \textbf{I} + \textbf{E}^{ij} + \textbf{E}^{ji}\right) \hat{\bo{\Xi}}_m^T,
\end{align}
and
\begin{align}\label{eq:mjade_kumulant2}
\hat{\textbf{C}}{}_{2, m}^{ij} := \hat{\textbf{B}}{}_{2, m}^{ij} - \hat{\bo{\Xi}}_m \left(\delta_{ij} \textbf{I} + \rho_m \textbf{E}^{ij} + \textbf{E}^{ji}\right) \hat{\bo{\Xi}}_m^T,
\end{align}
where $i,j = 1,\ldots,p_m$, and the estimates $\hat{\textbf{B}}{}_{1, m}^{ij}$, $\hat{\textbf{B}}{}_{2, m}^{ij}$ and $\hat{\bo{\Xi}}_m$ are obtained by applying the definitions of $\textbf{B}_{1, m}^{ij}$, $\textbf{B}_{2, m}^{ij}$ and $\bo{\Xi}_m$ to the empirical distribution of $\textbf{X}$, including an empirical standardization by $\hat{\bo{\Sigma}}{}_1^{-1/2},\ldots,\hat{\bo{\Sigma}}{}_r^{-1/2}$. Choosing then either of the sets, $c=1,2$, the rotation matrix $\textbf{U}_m^T$, $m=1,\ldots,r$ is found by simultaneous (approximate) diagonalization as
\begin{align}\label{eq:mjade_rotation}
\textbf{U}_m^T = \underset{\textbf{U}: \ \textbf{U}^T\textbf{U} = \textbf{I}}{\mbox{argmax}} \sum_{i=1}^{p_m} \sum_{j=1}^{p_m} \| \mbox{diag}(\textbf{U} \textbf{C}_{c, m}^{ij}(\textbf{X}_{st}) \textbf{U}^T) \|^2.
\end{align}
The corresponding estimates $\hat{\textbf{U}}{}_m^T$, $m=1,\ldots,r$, are obtained by replacing in \eqref{eq:mjade_rotation} the matrices $\textbf{C}_{c, m}^{ij}(\textbf{X}_{st})$ with their estimates $\hat{\textbf{C}}{}_{c, m}^{ij}$ .


Combining the standardization and the rotation, the final TJADE algorithm for a sample, $\textbf{X}_i \in \mathbb{R}^{p_1 \times \cdots \times p_r}$, $i=1,\ldots,n$, consists of the following steps.
\begin{enumerate}
\item[1)] Center $\textbf{X}_i$ and estimate $\hat{\bo{\Sigma}}_1,\ldots,\hat{\bo{\Sigma}}_r$.
\item[2)] Standardize: $\textbf{X}_i \leftarrow \textbf{X}_i \odot_1 \hat{\bo{\Sigma}}{}_1^{-1/2} \cdots \odot_r \hat{\bo{\Sigma}}{}_r^{-1/2}$.
\item[3)] Choose $c$ and estimate the $r$ rotations $\hat{\textbf{U}}{}_1^T,\ldots,\hat{\textbf{U}}{}_r^T$ by diagonalizing for each $m = 1,\ldots,r$ simultaneously the sets $\hat{\textbf{C}}{}_{c, m}^{ij}$, $i,j=1,\ldots,p_m$.
\item[4)] Rotate: $\textbf{X}_i \leftarrow \textbf{X}_i \odot_1 \hat{\textbf{U}}{}_1^T \cdots \odot_r \hat{\textbf{U}}{}_r^T$.
\end{enumerate}
Using Lemma 5.1.1 from \cite{virta2015mfobi} the final result can be written as the product $\textbf{X}_i \odot_1 \hat{\bo{\Phi}}_1 \cdots \odot_r \hat{\bo{\Phi}}_r$, where $\hat{\bo{\Phi}}_m := \hat{\textbf{U}}{}_m^T \hat{\boldsymbol{\Sigma}}{}_m^{-1/2}$, $m = 1,\ldots,r$, is the \textit{$m$-mode TJADE estimate}.
\begin{remark}
Technically, there is no reason why we could not use different $c$ for estimating different rotations $\textbf{U}_m$. However, the asymptotic properties of the different approaches are in the next section shown to be equivalent and thus the choice of $c$ is for large enough samples irrelevant.
\end{remark}

For a vector valued $\textbf{x}\in \mathbb{R}^p$ and a full-rank matrix $\textbf{A}\in \mathbb{R}^{p\times p}$, $(\textbf{A}\textbf{x})_{st}=\textbf{U}\textbf{x}_{st}$ for some orthogonal $\textbf{U}$ \cite{ilmonen2012invariant}. Unfortunately, the analogous relation in the tensor setting,
\begin{equation}\label{eq:tensor_ae}
(\textbf{X} \odot_1 \textbf{A}_1 \cdots \odot_r \textbf{A}_r)_{st} = \textbf{X}_{st} \odot_1 \textbf{U}_1 \cdots \odot_r \textbf{U}_r
\end{equation}
for some orthogonal $\textbf{U}_1,\ldots,\textbf{U}_r$, holds only for orthogonal $\textbf{A}_1,\ldots,\textbf{A}_r$. This lack of \textit{m-affine equivariance} of $\bo{\Sigma}_m(\textbf{X})$, $m=1,\ldots,r$, is discussed in \cite{virta2015mfobi} along with a conjecture that in the general tensor case, $r > 1$, no standardization functional exists which would lead into the property \eqref{eq:tensor_ae}. In practice this means that outside the model \eqref{eq:tica_model} a change (other than rotation or reflection) in the coordinate system leads into different estimated components. However, the TJADE estimator is still Fisher consistent by Theorem~\ref{theo:cumulant_diag}.

\section{Asymptotic properties}\label{sec:asymp}
The asymptotical properties of JADE were considered in \cite{BonhommeRobin2009}, \cite{miettinen2014fourth}, \cite{virta2015joint} and are in \cite{miettinen2014fourth}, \cite{virta2015joint} based on the fact that the JADE functional is affine equivariant, allowing them to consider only the case of no mixing, $\bo{\Omega} = \textbf{I}$. In the following we consider the analogous case of $\bo{\Omega}_1 = \textbf{I},\ldots,\bo{\Omega}_r = \textbf{I}$ for TJADE. However, because of the lack of full affine equivariance, the results generalize only to orthogonal mixing from all $r$ modes.

For a tensor $\textbf{X} \in \mathbb{R}^{p_1 \times \cdots \times p_r}$ define its \textit{$m$-flattening} $\textbf{X}_{(m)} \in \mathbb{R}^{p_m \times \rho_m}$ as the horizontal stacking of all $m$-mode vectors of the tensor into a matrix in a predefined order, see \cite{de2000multilinear} for a rigorous definition. If the stacking order is assumed to be cyclical in the dimensions in the sense of \cite{de2000multilinear} we have for $\textbf{X}^* := \textbf{X} \odot_1 \textbf{A}_1 \cdots \odot_r \textbf{A}_r$ the identity
\begin{equation}\label{eq:stacking_kronecker}
\textbf{X}^*_{(m)} = \textbf{A}_m \textbf{X}_{(m)} \left( \textbf{A}_{m+1} \otimes \cdots \otimes \textbf{A}_r \otimes \textbf{A}_1 \otimes \cdots \otimes \textbf{A}_{m-1} \right)^T.
\end{equation}
The reason why $m$-flattening is particularly useful for us is that it allows us to write the $m$-mode product of a tensor with itself as an ordinary matrix product, namely $\textbf{X} \odot_{-m} \textbf{X} = \textbf{X}_{(m)} \textbf{X}_{(m)}^T$, regardless of the stacking order. This, combined with the fact that the matrices $\textbf{C}{}^{ij}_{c, m}(\textbf{X})$ depend on $\textbf{X}$ only via the previous product, implies that it is sufficient to derive the asymptotics for the case $r = 2$ only. The results for tensors of order $r > 2$ are then obtained by applying the case $r = 2$ for each of the $m$-flattened matrices $\textbf{X}_{(1)},\ldots,\textbf{X}_{(r)}$. Similarly, even for the case $r = 2$ we only need to consider the $1$-mode TJADE estimate  $\hat{\bo{\Phi}}_1$ (matrix multiplication from left) as the results for $\hat{\bo{\Phi}}_2$ follow by simply transposing $\textbf{X}$ or, in the language of tensors, flattening $\textbf{X}$ from the second mode. Interestingly, we also have no need to specify the used set of cumulant matrices $c$, as the two choices, $c=1$ and $c=2$, are shown to lead into asymptotically equivalent estimators.

We next provide the asymptotic expressions for the elements of the TJADE estimate $\hat{\bo{\Phi}}_1 =: \hat{\bo{\Phi}}$ in the case of a matrix-valued sample $\textbf{X}_i \in \mathbb{R}^{p_1 \times p_2}$, $i=1, \ldots ,n$. 
The asymptotic properties of $\hat{\bo{\Phi}}$ can be shown to depend on row means of various moments of $\textbf{Z}$, particularly on the elements of $\bar{\bo{\kappa}}^{(1)}$ but also on the following
\begin{align*}
\bar{\bo{\beta}}^{(1)} &:= \frac{1}{p_2} \sum_{l=1}^{p_2} \left( \E[z_{1 l}^4], \ldots, \E[z_{p_1 l}^4] \right){}^T, \\
\bar{\bo{\omega}}^{(1)} &:= \frac{1}{p_2} \sum_{l=1}^{p_2} \left( \var[z_{1 l}^3], \ldots, \var[z_{p_1 l}^3] \right){}^T.
\end{align*}
Define further the covariance of two rows of kurtoses as
\[ \rho_{kk'} = \frac{1}{p_2}\sum_{l=1}^{p_2} \left(\beta_{kl} \beta_{k'l} \right)  - \bar{\beta}{}^{(1)}_k \bar{\beta}{}^{(1)}_{k'},\]
where $\beta_{kl} := \E[z_{kl}^4]$. To construct an asymptotic expression for $\hat{\bo{\Phi}}$ in Theorem~\ref{theo:mjade_asymp} we need the terms

\begin{align*}
\hat{s}_{kk'} &:= \frac{1}{p_2} \sum_{l=1}^{p_2} \left( \frac{1}{n} \sum_{i=1}^n z_{i,kl} z_{i,k'l} \right), \\
\hat{q}_{kk'} &:= \frac{1}{p_2} \sum_{l=1}^{p_2} \left( \frac{1}{n}\sum_{i=1}^n \left( z_{i,kl}^3 - E[z_{kl}^3] \right) z_{i,k'l} \right), \\
\hat{r}_{kk'} &:= \frac{1}{p_2} \sum_{l=1}^{p_2} \underset{l' \neq l}{\sum_{l'=1}^{p_2}} \left( \frac{1}{n}\sum_{i=1}^n z_{i,kl}^2 z_{i,kl'} z_{i,k'l'} \right),
\end{align*}
the joint limiting normality of which is easy to show, assuming the eighth moments of $\textbf{Z}$ exist.

\begin{theorem}\label{theo:mjade_asymp}
Let $\textbf{Z}_1,\ldots,\textbf{Z}_n$ be a random sample from a distribution with finite eighth moments and satisfying Assumptions \ref{assu:ind}, \ref{assu:stand} and \ref{assu:mjade_3} (see below). Then there exists a sequence of TJADE estimates such that $\hat{\boldsymbol{\Phi}} \rightarrow_P \textbf{\em I}$ and
\begin{align*}
\sqrt{n} (\hat{\phi}_{kk} - 1) &= -\dfrac{1}{2} \sqrt{n} (\hat{s}_{kk} - 1) + o_P(1), \\
\sqrt{n} \hat{\phi}_{kk'} &= \frac{\sqrt{n} \hat{\psi}_{kk'} + \sqrt{n} \hat{\psi}_{k'k} - d_{kk'} \sqrt{n} \hat{s}_{kk'} }{(\bar{\kappa}{}^{(1)}_k)^2 + (\bar{\kappa}{}^{(1)}_{k'})^2} + o_P(1),
\end{align*}
where $k \neq k'$, $\hat{\psi}_{kk'} := \bar{\kappa}{}^{(1)}_k (\hat{r}_{kk'} + \hat{q}_{kk'})$ and $d_{kk'} := (p_2 + 2)(\bar{\kappa}{}^{(1)}_k - \bar{\kappa}{}^{(1)}_{k'}) + (\bar{\kappa}{}^{(1)}_k)^2$.
\end{theorem}

Using the expressions of Theorem \ref{theo:mjade_asymp} the asymptotic variances of the elements of $\hat{\bo{\Phi}}$ can now be computed.

\begin{corollary}\label{cor:mjade_asymp}
Under the assumptions of Theorem \ref{theo:mjade_asymp} the limiting distribution of $\sqrt{n} \, \vec (\hat{\boldsymbol{\Phi}} - \textbf{\em I})$ is multivariate normal with mean vector $\textbf{0}$ and the following asymptotic variances.
\begin{alignat*}{3}
&ASV&(\hat{\phi}_{kk}) &= \frac{\bar{\beta}{}^{(1)}_k - 1}{4 p_2}, &\\
&ASV&(\hat{\phi}_{kk'}) &= \frac{\zeta_k + \zeta_{k'} + (\bar{\kappa}{}^{(1)}_{k'})^4 - 2 \bar{\kappa}{}^{(1)}_{k} \bar{\kappa}{}^{(1)}_{k'} \rho_{kk'}}{p_2 ((\bar{\kappa}{}^{(1)}_{k})^2 + (\bar{\kappa}{}^{(1)}_{k'})^2)^2}, &\quad k \neq k',
\end{alignat*}
where $\zeta_k := (\bar{\kappa}{}^{(1)}_{k})^2[\bar{\omega}{}^{(1)}_{k} - (\bar{\beta}{}^{(1)}_{k})^2] + (\bar{\kappa}{}^{(1)}_{k})^2 (\bar{\kappa}{}^{(1)}_{k} - 2)(p_2 - 1)$.
\end{corollary}

It is easily seen that the expressions in Corollary \ref{cor:mjade_asymp} revert to the forms of Corollary 4 in \cite{miettinen2014fourth} when $r=1$, that is, we observe just a vector $\textbf{x}$. In this case $\bar{\bo{\kappa}}^{(1)}$ contains just the element-wise kurtoses of the elements of $\textbf{z}$. Of the popular ICA methods, FastICA, FOBI and JADE, it is well-known that only for FOBI does the asymptotic behavior of $\hat{\phi}_{kk'}$ depend on components other than $z_k$ and $z_{k'}$. The analogous result holds also for TFOBI and TJADE in the sense that in TFOBI the asymptotic behavior of $\hat{\phi}^{(m)}_{kk'}$ depends on the whole tensor $\textbf{Z}$ \cite{virta2015mfobi} and in TJADE only on the $k$th and $k'$th $m$-mode faces of $\textbf{Z}$.

The denominators in Theorem \ref{theo:mjade_asymp} imply that for the existence of the limiting distributions we need the following assumption.

\begin{assumption}\label{assu:mjade_3}
For each $m = 1,\ldots,r$, at most one of the components of $\bar{\bo{\kappa}}^{(m)}$ is zero.
\end{assumption}

Assumption \ref{assu:mjade_3} for TJADE is much less restrictive than the assumption needed for TFOBI, for each $m=1,..,r$ the components of $\bar{\bo{\kappa}}^{(m)}$are distinct \cite{virta2015mfobi}, and the one needed for vector JADE, at most one element of $\bo{\kappa}$ is zero \cite{miettinen2014fourth}. More specifically, in TJADE, and in tensor independent component analysis in general, several individual elements of $\textbf{Z}$ are allowed to be Gaussian, as long as Assumption \ref{assu:gauss} is not violated. Conveniently located, a majority of the elements of $\textbf{Z}$ can thus be Gaussian.


The analytical comparison of TJADE and TFOBI via the asymptotic variances involves in general case rather complicated expressions and thus we resort to simulations for their comparison in the next section.

\section{Simulations and examples}\label{sec:simu}

In the following all computations were done in R 3.1.2 \cite{Rcore} especially using the R-packages JADE \cite{jade}, Rcpp \cite{eddelbuettel2011rcpp, eddelbuettel2013seamless} and ggplot2 \cite{ggplot2}. For the approximate joint diagonalization, an algorithm based on Jacobi rotations was used, see e.g \cite{belouchrani1997blind}. Testing the algorithms in various settings showed that both $c=1$ and $c=2$ yield almost identical results with respect to the MDI-values (see below) but the former is computationally more efficient and thus the TJADE solution in the simulations is computed with the choice $c=1$.

\subsection{Efficiency comparisons}

We compared the separation performance of TJADE with its nearest competitor, TFOBI, and also with regular FOBI and JADE as applied to vectorized tensor data, called here VFOBI and VJADE. Note that VFOBI and VJADE do not use the prior information on the data structure and are therefore expected to be worse than TFOBI and TJADE. The simulation setting was the same as in \cite{virta2015mfobi}: we simulated $n$ independent $3 \times 4$ matrix observations with individual elements coming from a diverse array of distributions. The excess kurtoses of the distributions used were -1.2, -0.6, 0, 1, 2, 3, 4, 5, 6, 8, 10 and 15 and the exact distributions used are given in Appendix \ref{sec:appen}.

We generated 2000 repetitions for each of the sample sizes, $n = 1000, 2000, 4000, 8000, 16000, 32000$, and for each sample the same data was mixed using three different distributions for the elements of the $1$-mode and $2$-mode mixing matrices, $\textbf{A}$ and $\textbf{B}$. In the first case the mixing matrices were random orthogonal matrices of sizes $3 \times 3$ and $4 \times 4$ distributed uniformly with respect to the Haar measure. In the second and third case the elements of both matrices were generated independently from $\mathcal{N}(0, 1)$ and $\mbox{Uniform}(-1, 1)$ distributions, respectively.

The mixed data were then subjected to each of the four different methods which produced the four unmixing matrix estimates, $\hat{\bo{\Phi}}_{VF}$, $(\hat{\bo{\Phi}}{}_{2,MF} \otimes \hat{\bo{\Phi}}{}_{1,MF})$, $\hat{\bo{\Phi}}_{VJ}$ and $(\hat{\bo{\Phi}}{}_{2,MJ} \otimes \hat{\bo{\Phi}}{}_{1,MJ})$. To allow comparing we took the Kronecker product of the $2$-mode and $1$-mode unmixing matrices of TFOBI and TJADE meaning that all the four previous matrices estimate the inverse of the same matrix $(\textbf{B} \otimes \textbf{A})$, up to scaling, sign-change and permutation of its columns.

The actual comparison was done by first computing the \textit{minimum distance index} (MDI) \cite{ilmonen2010new} of the estimates
\[D(\hat{\bo{\Phi}} \bo{\Omega}) = \frac{1}{\sqrt{p - 1}} \underset{\textbf{C} \in \mathcal{C}}{\text{inf}}\| \textbf{C} \hat{\bo{\Phi}} \bo{\Omega} - \textbf{I} \|_F, \]
where $\hat{\bo{\Phi}} \in \mathbb{R}^{p \times p}$ is the estimated unmixing matrix, $\bo{\Omega} \in \mathbb{R}^{p \times p}$ is the true mixing matrix and $\mathcal{C}$ is the set of all $p \times p$ matrices having a single non-zero element in each row and column. MDI thus measures how far away $\hat{\bo{\Phi}} \bo{\Omega}$ is from the set $\mathcal{C}$. The index varies from 0 to 1 with 0 indicating that the separation worked perfectly. In our simulation we further transformed the MDI-values as $n (p - 1) \mbox{MDI}^2$ which in vector-valued independent component analysis converges in distribution to a random variable with finite mean and variance \cite{ilmonen2010new}.

\begin{figure*}[!ht]
\centering
\includegraphics[width=\textwidth]{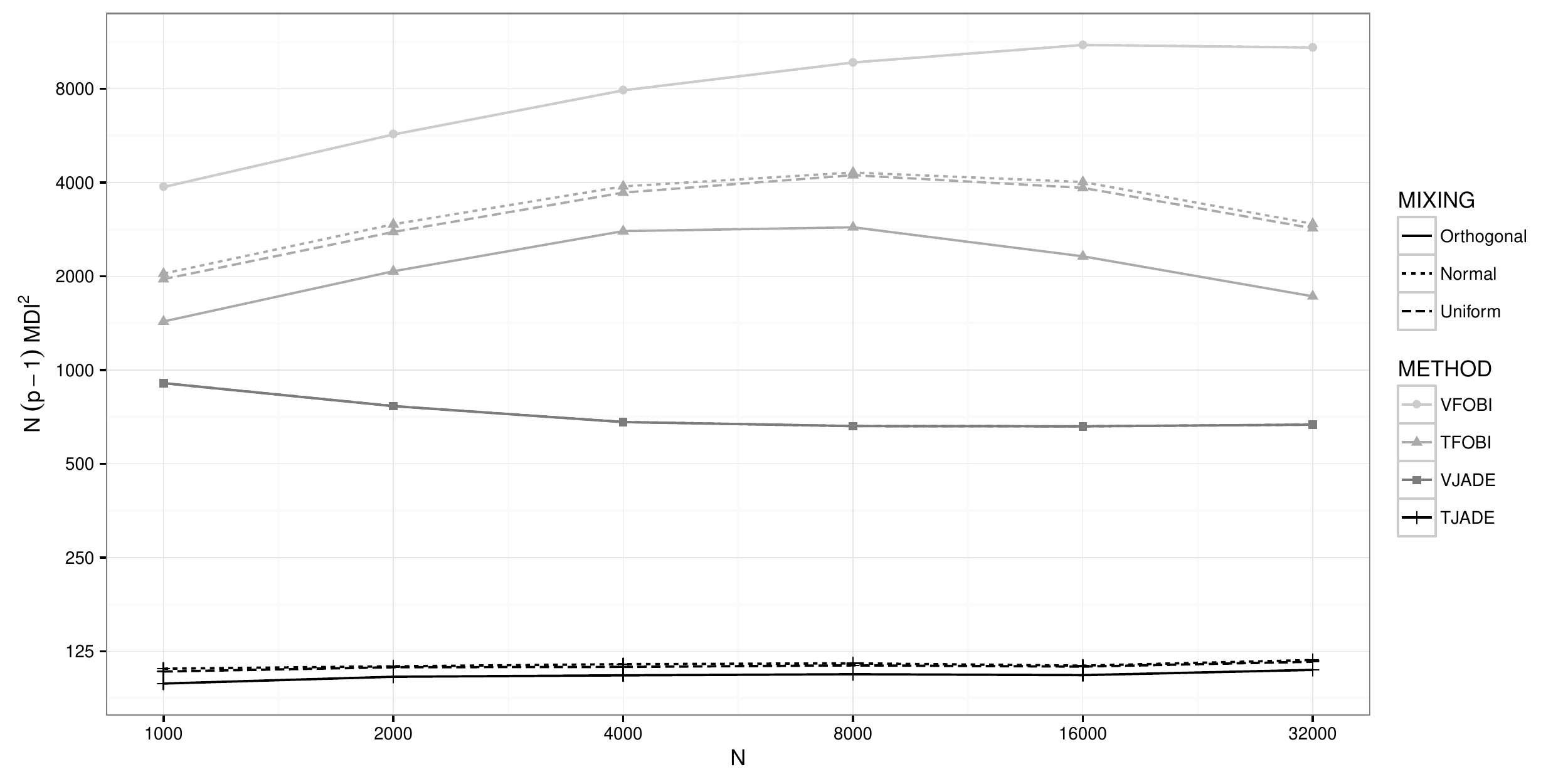}
\caption{Plot of sample size versus the transformation $n(p-1)$MDI$^2$ under combinations of the four different methods and three different distributions for the mixing matrices.}
\label{fig:mdi}
\end{figure*}

The mean transformed MDI-values for different sample sizes, methods and mixing matrices are shown in Figure \ref{fig:mdi}. The lines for both FOBI and JADE are for all mixings identical since both methods are affine equivariant. For TFOBI and TJADE the separation is best under orthogonal mixing, the results for normal and uniform mixing being a bit worse. But the main implication of the plot is that none of the other methods can really compete with TJADE in matrix independent component analysis. Interestingly, also regular JADE combined with vectorization is better than TFOBI.

\def\ttt{\texttt}

\subsection{Assumption comparisons}

As our second simulation study we compared the four methods of the previous simulation via their assumptions. For this we constructed three simulation settings of $3 \times 3 \times 2$ tensors with independent elements having either Gaussian ($\ttt N$), Laplace ($\ttt L$), exponential ($\ttt E$), or continuous uniform ($\ttt U$) distributions standardized to have zero means and unit variances. The distributions of the tensors are shown in the following by the two $3 \times 3 \times 1$ faces of each setting:
\begin{align*}
\mbox{Setting 1}: \quad &
\begin{pmatrix}
\texttt{N} & \texttt{L} & \texttt{E} \\
\texttt{L} & \texttt{L} & \texttt{E} \\
\texttt{E} & \texttt{E} & \texttt{E}
\end{pmatrix}
&
\begin{pmatrix}
\texttt{U} & \texttt{U} & \texttt{U} \\
\texttt{U} & \texttt{L} & \texttt{L} \\
\texttt{U} & \texttt{L} & \texttt{E}
\end{pmatrix} \\
\mbox{Setting 2}: \quad &
\begin{pmatrix}
\texttt{N} & \texttt{L} & \texttt{L} \\
\texttt{L} & \texttt{L} & \texttt{L} \\
\texttt{L} & \texttt{L} & \texttt{L}
\end{pmatrix}
&
\begin{pmatrix}
\texttt{U} & \texttt{U} & \texttt{U} \\
\texttt{U} & \texttt{L} & \texttt{L} \\
\texttt{U} & \texttt{L} & \texttt{L}
\end{pmatrix} \\
\mbox{Setting 3}: \quad &
\begin{pmatrix}
\texttt{E} & \texttt{E} & \texttt{N} \\
\texttt{E} & \texttt{E} & \texttt{N} \\
\texttt{N} & \texttt{N} & \texttt{N}
\end{pmatrix}
&
\begin{pmatrix}
\texttt{N} & \texttt{N} & \texttt{N} \\
\texttt{N} & \texttt{N} & \texttt{N} \\
\texttt{N} & \texttt{N} & \texttt{N}
\end{pmatrix}
\end{align*}
It is easy to see that none of the above settings satisfies the assumptions of VFOBI as all of them have at least two identical components. Only setting 1 satifies the assumption of TFOBI on distinct kurtosis means in all modes and settings 1 and 2 satisfy the assumption on maximally one component having zero excess kurtosis required by VJADE. All three settings satisfy Assumption \ref{assu:mjade_3} on maximally one zero kurtosis mean in each mode required by TJADE.

\begin{figure*}[!t]
\centering
\includegraphics[width=\textwidth]{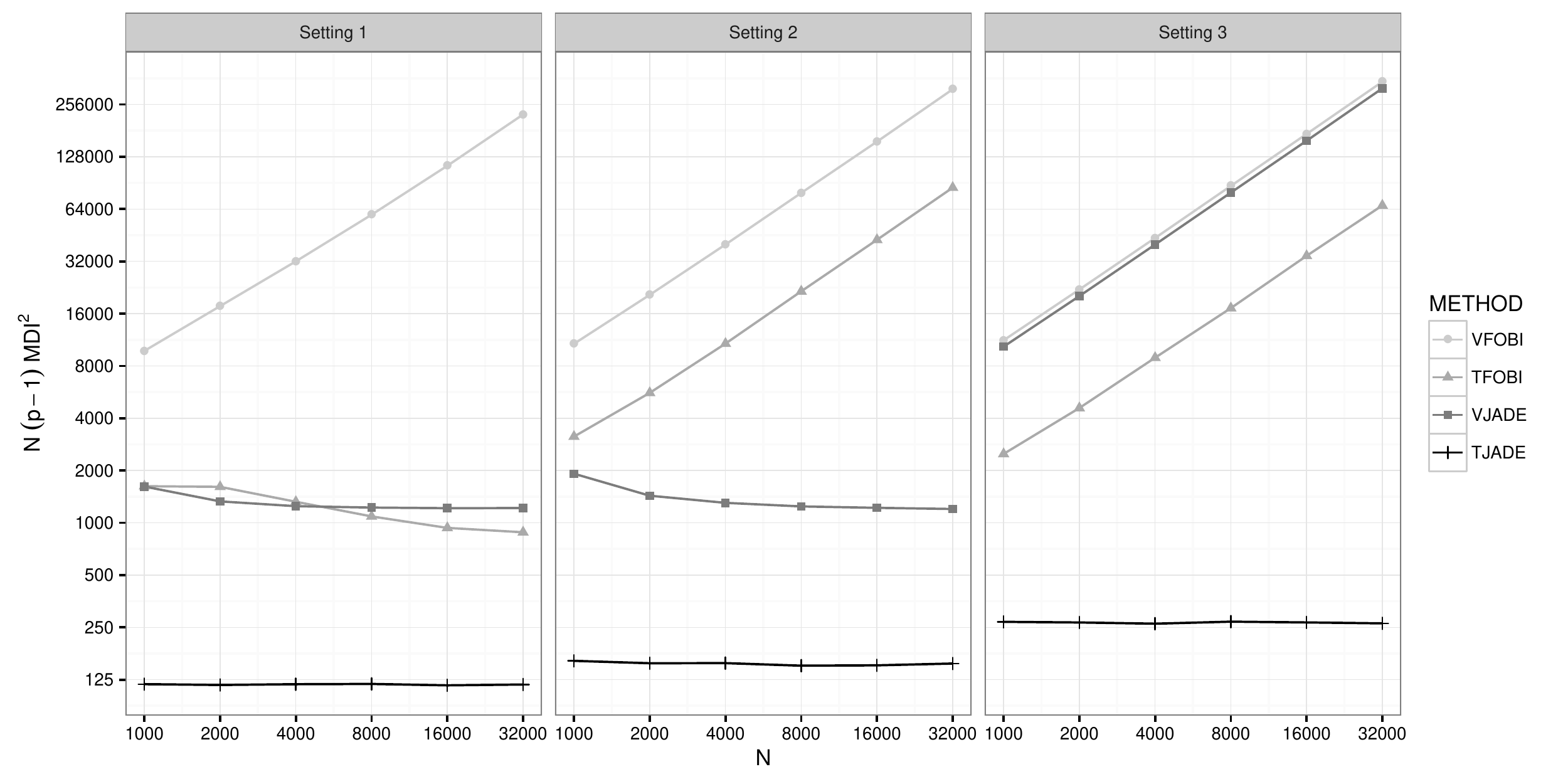}
\caption{Means of transformed MDI-values for different combinations of setting, sample-size and method. Moving from left to right, all other methods but TJADE break down one-by-one.}
\label{fig:assu}
\end{figure*}

We simulated 2000 repetitions of all three settings for different sample sizes using identity mixing and the resulting transformed MDI-values of the four methods are depicted in Figure \ref{fig:assu}. The above reasoning about the violation of assumptions is clearly visible in the plots. The mean transformed MDI-values of the different methods break one-by-one when the setting changes from 1 to 2 to 3 leaving TJADE as the only method able to handle all three settings. Interestingly, VJADE failed to converge 4601 times out of the 36000 total repetitions across all settings, the majority of failures occurring in the third setting.

The plot for setting 1 further indicates that there exist cases where TFOBI beats VJADE, proving that, though very efficient, the JADE methodology itself is not the only factor in the superior performance of TJADE; the tensor structure also plays an important role.

\subsection{Real data example}

Extreme kurtosis can be shown to be associated with multimodal distributions and thus independent component analysis is commonly used as preprocessing step in classification to obtain directions of interest. In this spirit we consider the \textit{semeion}\footnote{Semeion Research Center of Sciences of Communication, via Sersale 117, 00128 Rome, Italy; Tattile Via Gaetano Donizetti, 1-3-5,25030 Mairano (Brescia), Italy. } data set, available in the UCI Machine Learning Repository \cite{Lichman:2013} as a classification problem. The data consist of 1593 binary $16 \times 16$ pixel images of hand-written digits. For this example we chose only the images representing the digits 0, 1 and 7, having respective group sizes of 161, 162 and 158. The objective is to find a few components separating the three digits.

\begin{figure*}[!t]
\centering
\includegraphics[width=\textwidth]{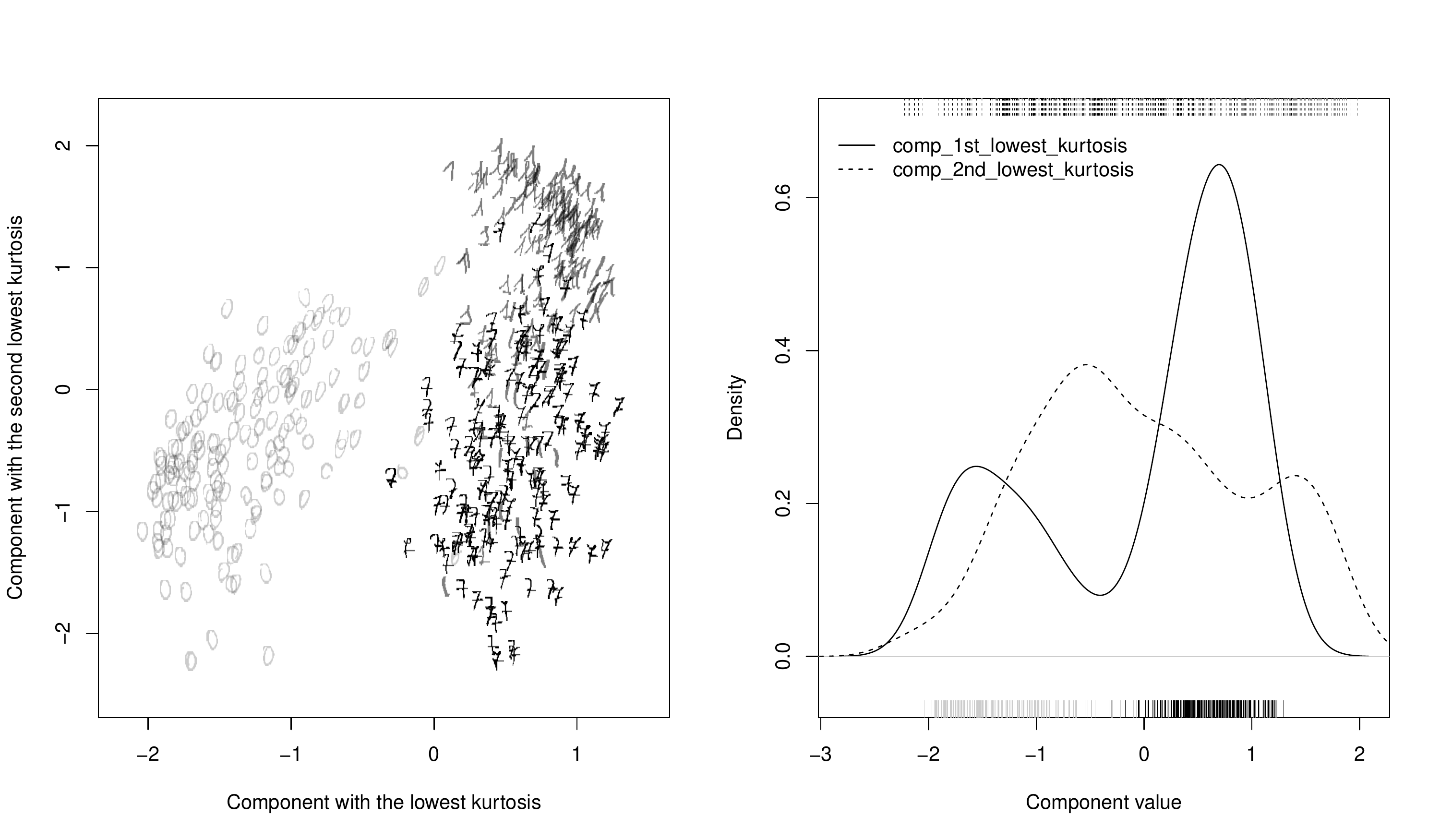}
\caption{Results of applying TJADE on the semeion data. The plot on the left-hand side shows the scatter plot of the two components having the lowest kurtoses found by TJADE with the individual images as markers. The three digits clearly form three groups in the plane. The density plots along with the rugs on the right-hand side imply the same. The lower rug corresponds to the component with the lowest kurtosis (min$\_$kurtosis$\_$1) and the coloring of the groups in the rugs is the same as in the scatter plot.}
\label{fig:semeion}
\end{figure*}

Subjecting the data to TJADE gives the results depicted in Figure \ref{fig:semeion}. The left-hand side plot shows the scatter plot of the two resulting components with the lowest kurtoses using the individual digit images as plot markers. Clearly the two found directions are sufficient to separate all three groups of digits. The same conclusion can be drawn from the corresponding density estimators and rug plots on the right-hand side of Figure \ref{fig:semeion}. As a next step, some low-dimensional classification algorithm could be applied to the extracted components to create a classification rule.

\section{Discussion}\label{sec:disc}

In this paper we proposed TJADE, an extension of the classic JADE suited for tensor-valued observations. Based on the same idea of diagonalizing multiple cumulant matrices as JADE, TJADE was shown to be very effective in solving the independent component problem. In the course of the paper we first reviewed the theory and the algorithm behind JADE and then formulated TJADE analogously giving two different, although asymptotically equivalent, ways of estimating the needed rotations. The asymptotic behaviors of the elements of the TJADE-estimates under orthogonal mixing were next provided allowing theoretical comparison to other methods. Finally, simulation studies comparing TJADE to TFOBI, and the na{\"i}ve approaches combining vectorization with either FOBI or JADE showed that TJADE is superior to all the previous competitors in tensor independent component analysis.

As the investigation of ICA methods for non-vector-valued objects is still in an early stage,  much further research is needed. Below we outline some ideas planned to follow this work.  

As the number of matrices to jointly diagonalize in estimating the $m$-mode rotation in TJADE grows proportional to the square of the corresponding dimension $p_m$, an extension like $k$-JADE \cite{MiettinenNordhausenOjaTaskinen:2013} is worth considering also for TJADE. And as a competing alternative also a tensor version of the FastICA algorithm \cite{HyvarinenKarhunenOja:2001} will be investigated. This opens many possibilities allowing choosing both the non-linearity function $g$ and the norm used in the maximization problem, see \cite{miettinen2015squared}.

\appendix[]\label{sec:appen}

The distributions used in the first simulation of Section \ref{sec:simu} are, starting from the upper left corner of the matrix and moving down and right, $\mbox{Uniform}(-\sqrt{3}, \sqrt{3})$, $\mbox{Triangular}(-\sqrt{6}, \sqrt{6}, 0)$, $\mathcal{N}(0, 1)$, $t_{10}$, $\mbox{Gamma}(3, \sqrt{3})$, $\mbox{Laplace}(0, 1/\sqrt{2})$, $\chi_3^2$, $\mbox{Gamma}(1.2, \sqrt{1.2})$, $\mbox{Exp}(1)$, $\chi_{1.5}^2$, $\chi_{1.2}^2$ and $\mbox{InverseGaussian}(1, 1)$. The distributions were further shifted to have zero means and unit variances.

\begin{proof}[The proof of Theorem \ref{theo:cumulant_diag}]

Consider first the case $c=1$ and the four terms in \eqref{eq:mjade_cij1} separately fixing the choice of $m$. Denoting the first term of \eqref{eq:mjade_cij1} by $\textbf{B}_{1, m}^{ij}(\textbf{X})$, then according to Lemma 5.4.1 in \cite{virta2015mfobi} we have
\begin{align*}
\ali \textbf{B}_{1, m}^{ij}(\textbf{X}_{st}) \\
\ali= \frac{\tau^4}{\rho_m} \textbf{U}_m E\left[ (\textbf{u}^{(m)}_i)^T \textbf{Z}_{(m)} \textbf{Z}{}_{(m)}^T \textbf{u}^{(m)}_j \cdot \textbf{Z}_{(m)} \textbf{Z}{}_{(m)}^T \right] \textbf{U}{}_m^T,
\end{align*}
where $\textbf{Z}_{(m)}$ is the flattened matrix defined in Section \ref{sec:asymp} and $(\textbf{u}^{(m)}_i)^T$ is the $i$th row of $\textbf{U}_m$. Using the standard properties of expected value and independent random variables the $(k,k')$ element of the inner expectation can be shown to be for $k \neq k'$ equal to $u^{(m)}_{ik} u^{(m)}_{jk'} + u^{(m)}_{jk} u^{(m)}_{ik'}$ and for $k = k'$ equal to $\delta_{ij} \rho_m + u^{(m)}_{ik} u^{(m)}_{jk} (\bar{\kappa}{}^{(m)}_k + 2)$. Using these to construct a matrix form for the expectation we have
\begin{align*}
\textbf{B}_{1, m}^{ij}(\textbf{X}_{st}) &= \tau^4 \textbf{U}_m \left( \sum_{k=1}^p u{}^{(m)}_{ik} u{}^{(m)}_{jk} \bar{\kappa}{}^{(m)}_k \textbf{E}^{kk} \right) \textbf{U}{}_m^T \\
& + \tau^4 \delta_{ij} \rho_m \textbf{I} + \tau^4 \textbf{E}^{ij} + \tau^4 \textbf{E}^{ji}.
\end{align*}
The second, third and fourth terms in \eqref{eq:mjade_cij1} then serve to remove the extra constant terms above. That they indeed cancel one-by-one the final terms can easily be shown by examining them in the above manner using the independence of $\textbf{X}$ and $\textbf{X}^*$. This concludes the proof for $c=1$ and the corresponding result for $c=2$ can be proven in precisely the same manner.
\end{proof}

\begin{proof}[The proof of Theorem \ref{theo:mjade_asymp}]
The consistency of the TJADE estimator is proven similarly as the consistency of the TFOBI estimator in the proof of Theorem 5.2.1 in \cite{virta2015mfobi}.

In the following we assume that $r = 2$ and we are interested in the asymptotical behavior of the $1$-mode unmixing matrix. As discussed in Section \ref{sec:asymp}, for the general case of arbitrary $r$ and $m$-mode unmixing matrix, it suffices to $m$-flatten the tensor and replace in the following $\hat{\bo{\Sigma}}{}_1^{-1/2}$ with $\hat{\bo{\Sigma}}{}_m^{-1/2}$, $\hat{\bo{\Sigma}}{}_2^{-1/2}$ with $\hat{\bo{\Sigma}}{}_{m+1}^{-1/2} \otimes \cdots \otimes \hat{\bo{\Sigma}}{}_{r}^{-1/2} \otimes \hat{\bo{\Sigma}}{}_{1}^{-1/2} \otimes \cdots \otimes \hat{\bo{\Sigma}}{}_{m-1}^{-1/2}$, $p_2$ with $\rho_m$ and use the corresponding row mean quantities.

For the asymptotic expressions of the diagonal elements of $\sqrt{n}(\hat{\bo{\Phi}} - \textbf{I})$ it suffices to use the same arguments as in the proof of Theorem 5.2.1 in \cite{virta2015mfobi} and for the off-diagonal elements we aim to use Lemma 2 from \cite{miettinen2014fourth}.

But first, define the \textit{symmetric} standardization functionals $\hat{\textbf{L}} = (\hat{l}_{kk'}) := \hat{\bo{\Sigma}}{}_1^{-1/2}$ and $\hat{\textbf{R}} = (\hat{r}_{ll'}) := \hat{\bo{\Sigma}}{}_2^{-1/2}$ giving the standardized identity-mixed observations as $\textbf{X}{}_{st, i} = \hat{\textbf{L}} \tilde{\textbf{Z}}_i \hat{\textbf{R}}{}^T$, where $\tilde{\textbf{Z}}_i = \textbf{Z}_i - \bar{\textbf{Z}} $. We then have
\[\sqrt{n}(\hat{l}_{kk'} - \delta_{kk'}) = -(1/2) \sqrt{n} (\hat{s}{}_{kk'} - \delta_{kk'}) + o_P(1),\]
see \cite{virta2015mfobi}, and as simple moment-based estimators we have both $\sqrt{n}(\hat{\textbf{L}} - \textbf{I}) = O_P(1)$ and $\sqrt{n}(\hat{\textbf{R}} - \textbf{I}) = O_P(1)$, regardless of whether we really have $r = 2$ or use flattened tensors of higher order.

Assume then first that $c=1$. The matrices $\hat{\textbf{C}}{}_{1, 1}^{kk'}$, $k,k' = 1,\ldots,p$, in \eqref{eq:mjade_kumulant1} to be simultaneously diagonalized satisfy $\hat{\textbf{C}}{}^{kk'} := \hat{\textbf{C}}{}_{1, 1}^{kk'} \rightarrow_P \textbf{C}{}_{1, 1}^{kk'}(\textbf{Z}_i) = \delta_{kk'} \bar{\kappa}{}^{(1)}_k \textbf{E}^{kk}$. In the view of Lemma 2 in \cite{miettinen2014fourth} this means that the only matrices $\textbf{C}{}_{1, 1}^{rs}(\textbf{Z}_i)$, $r,s=1,\ldots,p$ having non-zero $k$th or $k'$th diagonal elements are $\textbf{C}{}_{1, 1}^{kk}(\textbf{Z}_i)$ and $\textbf{C}{}_{1, 1}^{k'k'}(\textbf{Z}_i)$, respectively, yielding the following form for the $(k, k')$, $k \neq k'$, element of the rotation $\hat{\textbf{U}} := \hat{\textbf{U}}{}_1^T$ estimated by \eqref{eq:mjade_rotation}.
\[
\sqrt{n}\hat{u}_{kk'} = \frac{\bar{\kappa}{}^{(1)}_k \sqrt{n} \hat{\textbf{C}}{}^{kk}_{kk'} - \bar{\kappa}{}^{(1)}_{k'} \sqrt{n} \hat{\textbf{C}}{}^{k'k'}_{kk'}}{(\bar{\kappa}{}^{(1)}_k)^2 + (\bar{\kappa}{}^{(1)}_{k'})^2} + o_P(1),
\]
where $\hat{\textbf{C}}{}^{kk}_{rs}$ is the $(r,s)$ element of $\hat{\textbf{C}}{}^{kk}$. The above expression then together with the $(k, k')$, $k \neq k'$, element of the left standardization matrix $\hat{\textbf{L}}$ gives an asymptotic expression for the off-diagonal elements of the estimated left TJADE matrix, see \cite{virta2015mfobi}:
\begin{align}\label{eq:proof_twopart}
\sqrt{n}\hat{\phi}_{kk'} = \sqrt{n}\hat{u}_{kk'} + \sqrt{n} \hat{l}_{kk'} + o_P(1),
\end{align}
reducing the problem of finding the asymptotics of TJADE into the task of finding the asymptotic behaviors of $\sqrt{n} \hat{\textbf{C}}{}^{kk}_{kk'}$ and $\sqrt{n} \hat{\textbf{C}}{}^{k'k'}_{kk'}$. Dropping the subscripts for clarity, note that $\hat{\textbf{C}}{}^{aa} = \hat{\textbf{B}}{}^{aa} - \hat{\bo{\Xi}} (p_2 \textbf{I} + 2 \textbf{E}^{aa}) \hat{\bo{\Xi}}{}^T$ and starting from $ \hat{\textbf{B}}{}^{aa}$ write it out as
\[\hat{\textbf{B}}{}^{aa} = \frac{1}{p_2 n}\sum_{i=1}^n (\hat{\textbf{L}}{}_a^T \tilde{\textbf{Z}}_i \hat{\textbf{R}}{}^* \tilde{\textbf{Z}}{}_i^T \hat{\textbf{L}}{}_a) \cdot \hat{\textbf{L}} \tilde{\textbf{Z}}_i \hat{\textbf{R}}{}^* \tilde{\textbf{Z}}{}_i^T \hat{\textbf{L}}{}^T,\]
where $\hat{\textbf{L}}{}_a^T$ is the $a$th row of $\hat{\textbf{L}}$ and $\hat{\textbf{R}}{}^* := \hat{\textbf{R}}{}^T\hat{\textbf{R}}$. An arbitrary off-diagonal element of $\sqrt{n}(\hat{\textbf{B}}{}^{aa} - \textbf{B}{}^{aa}(\textbf{Z}_i))$ then has after the matrix multiplication the form
\begin{align}\label{eq:proof_asyexp}
\sqrt{n}\hat{\textbf{B}}{}^{aa}_{kk'} = \frac{1}{p_2 n}\sum_{defgstuv} \sqrt{n} \hat{r}^*_{ef} \hat{r}^*_{tu}  \hat{l}_{ad} \hat{l}_{ag} \hat{l}_{ks}  \hat{l}_{k'v} \hat{H}_{de, gf, st, vu},
\end{align}
where $\hat{H}_{de, gf, st, vu} = (1/n) \sum_{i=1}^n \tilde{z}_{i, de} \tilde{z}_{i, gf} \tilde{z}_{i, st} \tilde{z}_{i, vu} \rightarrow_P E(z_{i, de} z_{i, gf} z_{i, st} z_{i, vu})$. Next we expand the multiplicands $\hat{r}{}^*_{\cdot\cdot}$ and $\hat{l}_{\cdot\cdot}$ in \eqref{eq:proof_asyexp} one-by-one such as $\hat{l}_{ab} = (\hat{l}_{ab} - \delta_{ab}) + \delta_{ab}$, the first term of which is $O_P(1)$ when combined with $\sqrt{n}$ allowing the use of Slutsky's theorem  to the whole multiple sum and the second term of which produces an expression like \eqref{eq:proof_asyexp} only with one summation index less.

Starting from left this process then produces the terms $o_P(1)$; $o_P(1)$; $\delta_{ak} \sqrt{n} \hat{l}_{kk'} + \delta_{ak'} \sqrt{n} \hat{l}_{k'k} + o_P(1)$; $\delta_{ak} \sqrt{n} \hat{l}_{kk'} + \delta_{ak'} \sqrt{n} \hat{l}_{k'k}  + o_P(1)$; $\delta_{ak'}(\bar{\kappa}{}^{(1)}_{k'} + p_2 + 2)\sqrt{n}\hat{l}_{kk'} + (1-\delta_{ak'}) p_2 \sqrt{n} \hat{l}_{kk'}  + o_P(1)$ and $\delta_{ak}(\bar{\kappa}{}^{(1)}_{k} + p_2 + 2)\sqrt{n}\hat{l}_{k'k} + (1-\delta_{ak}) p_2 \sqrt{n} \hat{l}_{k'k} + o_P(1)$ finally leaving us with the expression
\begin{align}\label{eq:proof_finalform}
\frac{1}{p_2}\sum_{et} \frac{1}{\sqrt{n}} \sum_{i=1}^n \tilde{z}{}_{i, ae}^2 \tilde{z}_{i, kt} \tilde{z}_{i, k't} + o_P(1).
\end{align}
Substituting now either $a=k$ or $a=k'$, expanding $\tilde{z}_{i, ab} = z_{i, ab} - \bar{z}_{ab}$ and using the quantities defined in Section \ref{sec:asymp} the expression in \eqref{eq:proof_finalform} gets the forms $\sqrt{n} \hat{r}_{kk'} + \sqrt{n} \hat{q}_{kk'} + o_P(1)$ and $\sqrt{n} \hat{r}_{k'k} + \sqrt{n} \hat{q}_{k'k} + o_P(1)$, respectively.

Using the above, e.g. $\sqrt{n}\hat{\textbf{B}}{}^{kk}_{kk'}$ gets the form
\[(p_2+2)\sqrt{n}\hat{l}_{kk'} + (\bar{\kappa}{}^{(1)}_{k} + p_2 + 2) \sqrt{n} \hat{l}_{k'k} + \sqrt{n} \hat{r}_{kk'} + \sqrt{n} \hat{q}_{kk'} + o_P(1).\]

For the asymptotic behavior of the remaining term $\hat{\bo{\Xi}} (p_2 \textbf{I} + 2 \textbf{E}^{aa}) \hat{\bo{\Xi}}{}^T$ one can first use techniques similar to the above to show for $\hat{\bo{\Xi}} = (\hat{\xi}_{kk'})$ that $\sqrt{n}(\hat{\xi}_{kk'} - \delta_{kk'}) = o_P(1)$ for $k \neq k'$. Consequently an arbitrary off-diagonal element of $\sqrt{n}(\hat{\bo{\Xi}} (p_2 \textbf{I} + 2 \textbf{E}^{aa}) \hat{\bo{\Xi}}{}^T - p_2 \textbf{I} - 2 \textbf{E}^{aa})$ is also $o_P(1)$ implying that the term actually contributes nothing to the asymptotic variances of the estimator. Thus $\sqrt{n}\hat{\textbf{C}}{}^{aa}_{kk'} = \sqrt{n}\hat{\textbf{B}}{}^{aa}_{kk'} + o_P(1)$ and the result of Theorem \ref{theo:mjade_asymp} is obtained by plugging everything in into \eqref{eq:proof_twopart} and using the fact that the standardization functionals are symmetric. The asymptotic variances of Corollary \ref{cor:mjade_asymp} are then straightforward to obtain, e.g. using the table of covariances in the proof of Theorem 5.2.1 in \cite{virta2015mfobi}.

Although the starting expressions for $c=1$ and $c=2$ are different the final expressions for both $\sqrt{n}\hat{\textbf{C}}{}^{kk}_{kk'}$ and $\sqrt{n}\hat{\textbf{C}}{}^{k'k'}_{kk'}$ actually match exactly. The corresponding proof for $c=2$ is then obtained in exactly likewise manner, expanding the terms suitably and using Slustky's theorem, and is thus omitted here.
\end{proof}

\section*{Acknowledgments}
The work of Virta, Nordhausen and Oja was supported by the Academy of Finland Grant 268703. The work of Li was supported by the National Science Foundation Grant DMS-1407537.

\ifCLASSOPTIONcaptionsoff
  \newpage
\fi



\bibliographystyle{IEEEtran}
\bibliography{new_references}
\end{document}